\documentclass[letterpaper, 10 pt, conference]{ieeeconf}  

\usepackage{amsmath,amssymb,amsfonts}
\usepackage{algorithm,algorithmic}
\usepackage{graphicx}
\usepackage{textcomp}
\usepackage{url}
\usepackage{color}

\usepackage{lipsum}
\usepackage{epstopdf}

\usepackage{amsopn}

\usepackage{amsthm}

\usepackage{enumerate}
\usepackage{cancel}
\usepackage{mathrsfs}
\usepackage{mathdots}
\usepackage{euscript}
\usepackage{amscd}
\usepackage{placeins}

\usepackage{tikz}
\usetikzlibrary{arrows}
\usepackage{verbatim}

\newtheorem{theorem}{Theorem}

\newtheorem{lemma}{Lemma}
\newtheorem{proposition}{Proposition}
\newtheorem{assumption}{Assumption}

\theoremstyle{definition}

\newtheorem{example}{Example}
\newtheorem{problem}{Problem}

\theoremstyle{remark}
\newtheorem{remark}{Remark}

\newcommand{\nn}{\nonumber}

\newcommand{\bmat}{\left[ \begin{matrix}}
	\newcommand{\emat}{\end{matrix} \right]}
\newcommand{\innerprod}[2]{\langle{#1},\,{#2}\rangle}

\DeclareMathOperator{\E}{{\mathbb E}}
\newcommand{\Rbb}{\mathbb R}
\newcommand{\Hbb}{\mathbb H}

\newcommand{\Zbb}{\mathbb Z}

\newcommand{\Tbb}{\mathbb T}

\newcommand{\zb}{\mathbf  z}

\newcommand{\ab}{\mathbf a}
\newcommand{\bb}{\mathbf  b}

\newcommand{\pb}{\mathbf  p}

\newcommand{\mb}{\mathbf  m}

\newcommand{\cb}{\mathbf c}
\newcommand{\qb}{\mathbf q}
\newcommand{\tb}{\mathbf t}
\newcommand{\kb}{\mathbf k}

\newcommand{\oneb}{\mathbf 1}
\newcommand{\zerob}{\mathbf 0}

\newcommand{\thetab}{\boldsymbol{\theta}}

\newcommand{\Pfrak}{\mathfrak{P}}

\renewcommand{\d}{\mathrm{d}}

\newcommand{\m}{\mu} 

\IEEEoverridecommandlockouts                              
\title{\LARGE \bf  A Weaker Regularity Condition \\ for the Multidimensional $\nu$-Moment Problem
}

\author{Bin Zhu and Mattia Zorzi
\thanks{This work was supported in part by the National Natural Science Foundation of China with grant 62103453, and the ``Hundred Talent'' Program of Sun Yat-sen University.}
\thanks{B. Zhu is with School of Intelligent Systems Engineering, Sun Yat-sen University, Gongchang Road 66, 518107 Shenzhen, China {\tt\small zhub26@mail.sysu.edu.cn}}%
\thanks{M. Zorzi is with Department of Information Engineering, University of Padova, Via Gradenigo 6/B, 35131 Padova, Italy {\tt\small zorzimat@dei.unipd.it}}%
}

\begin{document}

\maketitle
\thispagestyle{empty}
\pagestyle{empty}

\begin{abstract}
We consider the problem of finding a $d$-dimensional spectral density through a moment problem which is characterized by an integer parameter $\nu$. Previous results showed that there exists an approximate solution under the regularity condition $\nu\geq d/2+1$. To realize the process corresponding to such a spectral density, one would take $\nu$ as small as possible.  In this letter we show that this condition can be weaken as $\nu\geq d/2$.
\end{abstract}

\section{Introduction}
Multidimensional stationary processes (or stationary random fields) represent a fundamental tool in many applications of signal and image processing. For instance, the data collected from an automotive radar system can be modeled by a multidimensional process of dimension $d=3$, see \cite{engels2014target,ZFKZ2019fusion}. In those applications we have to   estimate the multidimensional spectral density of the process. This task can be addressed by means of a moment problem, more precisely, a \emph{convex optimization problem subject to moment constraints}.

In the unidimensional case $(d=1)$ a wide range of spectral estimation paradigms based on moment problems have been proposed, see for instance \cite{BGL-98,FPR-08}: in this case the moments correspond to some covariance lags of the process and in the simplest setup the optimal spectrum maximizes the entropy rate. The appealing property of these paradigms is that the optimal spectrum is rational and thus leading to a finite-dimensional linear stochastic system (called ``shaping filter'' in the literature of signal processing) after spectral factorization. 

In the case where the moments include both covariance lags and cepstral coefficients (i.e., logarithmic moments), then it is possible to characterize only an approximate solution to the moment problem: the spectrum maximizing the entropy rate matches the covariance lags and approximately the cepstral coefficients  \cite{enqvist2004aconvex}. Such a solution is obtained by considering a \emph{regularized} version of the dual optimization problem.

These paradigms have been extended also to the multidimensional case, see e.g., \cite{RKL-16multidimensional,ringh2018multidimensional}. Although spectral factorization is not always possible in the multidimensional setting, rationality still seems to be a key ingredient toward a finite-dimensional realization theory \cite{geronimo2006factorization,geronimo2004positive}. The main issue is, however, that the solution of the moment problem is not necessarily a spectral density, but rather a \emph{spectral measure} that may contain a singular part \cite{KLR-16multidimensional}. In particular, if the moments are both the covariance lags and the cepstral coefficients, the existence of an approximate rational solution is only guaranteed when the dimension is $d\leq 2$ \cite{KLR-16multidimensional}.

In order to overcome this limitation on the dimension $d$, we have proposed a new moment problem, hereafter called \emph{$\nu$-moment problem}, in which the entropy rate has been replaced by a more general definition of entropy, called \emph{$\nu$-entropy}, whose derivation comes from the $\alpha$-divergence \cite{Z-14rat}. The definition of cepstral coefficients has been generalized accordingly.  
The $\nu$-moment problem is characterized by the integer parameter $\nu$. In \cite{Zhu-Zorzi2023cepstral,Zhu-Zorzi2022ceps_est} we have shown that for any $d>2$ there exists an approximate rational solution to the $\nu$-moment problem under the regularity condition $\nu\geq d/2+1$.
On the other hand, if the solution admits a spectral factorization, then the estimated process can be realized through a cascade of $\nu$ identical linear filters. Therefore, the larger $\nu$ is, the larger the complexity is in order to realize such a process. Accordingly, from a practical perspective one would take  $\nu$ as small as possible.

The aim of the present letter is to show that the regularity condition can be weakened as $\nu\geq d/2$ and thus it is possible to estimate a $d$-dimensional process with $\nu=d/2$ whose realization, if admissible, is simpler than the one obtained using the theory in \cite{Zhu-Zorzi2023cepstral}. Such a result is achieved through a new regularization technique. Moreover, the technical proof of the existence of such a solution takes a different route from the one in \cite{Zhu-Zorzi2023cepstral}.

The outline of this letter is as follows. In Section \ref{sec:intr} we introduce the multidimensional $\nu$-moment problem and the regularization term for the dual problem, while in Section \ref{sec:dual} we characterize the regularized dual problem. In Section
\ref{sec:weaker_cond} we prove the existence of an approximate solution for the $\nu$-moment problem for $\nu\geq d/2$. Some numerical examples are provided in Section \ref{sec:nume}. Finally, in Section \ref{sec:concl} we draw the conclusions.

\section{Problem formulation} \label{sec:intr}

Consider a $d$-dimensional real\footnote{We choose to present the theory for real random fields for simplicity. After suitable adaptation, complex random fields can be handled as well.} stationary random field $\{y(\tb) : \tb\in\Zbb^d\}$ with zero mean and a spectral density $\Phi(\thetab)$. The latter is a nonnegative function on the $d$-dimensional frequency domain $\Tbb^d:=(0,2\pi]^d$ and $\thetab=(\theta_1, \dots, \theta_d)\in\Tbb^d$ is a frequency vector.
In the sequel, we shall use the notation $\Phi(e^{i \thetab})$ which is common in Complex Analysis. Here $e^{i\thetab}$ is a shorthand for the vector $(e^{i\theta_1}, \dots, e^{i\theta_d})$ representing a point on the $d$-torus which is isomorphic to $\Tbb^d$.

The $\nu$-entropy of the random field $y$ with the integer parameter $\nu>1$ is defined as 
 \begin{align}\label{nu_entropy_formula}
 \Hbb_\nu(\Phi) := \frac{\nu^2}{\nu-1} \left(\int_{\Tbb^d}\Phi^{\frac{\nu-1}{\nu}}\d\m-1\right)
\end{align}
where 
$\d\mu=\frac{1}{(2\pi)^d}\prod_{j=1}^{d}\d\theta_j$
is the normalized Lebesgue measure on $\Tbb^d$. It is worth noting that $ \Hbb_\nu$ is a more general definition of entropy. Indeed, for the case $\nu= 1$ (which is understood in a suitable limit sense), we obtain the usual entropy rate $\Hbb_1(\Phi) = \int_{\Tbb^d}\log \Phi \d\m$ \cite{Zhu-Zorzi2023cepstral}. In this letter, we will consider the following multidimensional $\nu$-moment problem 
\begin{subequations}\label{primal_prob}
	\begin{alignat}{2}
		& \underset{ \Phi\geq 0}{\max}
		& \quad & \Hbb_\nu(\Phi)  \label{nu-entropy} \\
		&\hspace{0.4cm} \text{s.t.}
		& \quad & c_\kb=\int_{\Tbb^d}e^{i\innerprod{\kb}{\thetab}}\Phi \d\m \;\ \forall \kb\in\Lambda, \label{constraint_cov} \\
		& 
		& \quad & m_\kb= \frac{\nu}{\nu-1}\int_{\Tbb^d}e^{i\innerprod{\kb}{\thetab}}\Phi^{\frac{\nu-1}{\nu}} \d\m\;\ \forall \kb\in\Lambda_0. \label{constraint_cep}
	\end{alignat}
\end{subequations}
where $\innerprod{\kb}{\thetab}=\sum_{j=1}^{d} k_j\theta_j$ is the standard inner product in $\Rbb^d$ and the multidimensional exponential function is understood as $e^{i\innerprod{\kb}{\thetab}} = \prod_{j=1}^d e^{ik_j\theta_j}$. In the moment constraints we have:
 
\begin{itemize}
	\item $\cb=\{c_{\kb}\,:\, \kb\in\Lambda\}$ which is a covariance multisequence of the random field $y$, namely $c_{\kb}=\E[y(\tb+\kb) y(\tb)]$ where $\E$ denotes the expectation operator. The finite index set $\Lambda\subset\Zbb^d$ contains $\zerob$ and is symmetric with respect to the origin, i.e., $\kb\in\Lambda$ implies  $-\kb\in\Lambda$. 
	It is well known that the covariances are the Fourier coefficients of the spectral density $\Phi$ which is the meaning of \eqref{constraint_cov}.
	\item $\mb=\{m_{\kb}\, :\, \kb\in\Lambda_0\}$ which is a multisequence of generalized cepstral coefficients, called \emph{$\nu$-cepstral coefficients}, associated to the same random field $y$, see   \cite{Zhu-Zorzi2023cepstral}. It is a generalization of the classic logarithmic moments used in \cite{enqvist2004aconvex,RKL-16multidimensional}. Notice that, the notion of $\nu$-cepstral coefficients is consistent with the objective functional \eqref{nu_entropy_formula} employed in \eqref{primal_prob}. The index set $\Lambda_0:=\Lambda\backslash\{\zerob\}$ is such that $m_{\zerob}$ is excluded (for technical reasons). 
\end{itemize}

The $\nu$-moment problem \eqref{primal_prob}, which is also referred to as the primal optimization problem, can be used to perform spectral estimation. Assume that a dataset generated by $y$ has been collected. Then, it is possible to compute the sample estimates $\hat c_\kb$ and $\hat m_\kb$ of $c_\kb$ and $m_\kb$, respectively. Then, the estimate of $\Phi$ is the solution of \eqref{primal_prob} in which $c_\kb$ and $m_\kb$ are substituted by $\hat c_\kb$ and $\hat m_\kb$.

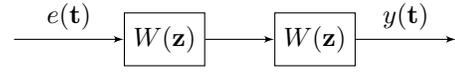
\begin{figure}[!t]
	\centering
	\tikzstyle{int}=[draw, minimum size=2em]
	\tikzstyle{init} = [pin edge={to-,thin,black}]
	\begin{tikzpicture}[node distance=2cm,auto,>=latex']
	\node [int] (a) {$W(\zb)$};
	\node (b) [left of=a, coordinate] {};
	\node [int] (c) [right of=a] {$W(\zb)$};
	\node (d) [right of=c] {};
	\path[->] (b) edge node {$e(\tb)$} (a);
	\path[->] (c) edge node {$y(\tb)$} (d);
	\draw[->] (a) edge node {} (c) ;
	\end{tikzpicture}
	\caption{A $d$-dimensional cascade linear stochastic system with two identical subsystems, where $\zb=(z_1, z_2, \dots, z_d)$.}
	\label{fig:cascade_linear_system}
\end{figure}

In \cite{Zhu-Zorzi2023cepstral}, we derived the dual optimization problem and showed that, \emph{if the dual problem admits an interior-point solution}, the optimal spectral density (primal variable) is a rational function of the form 
\begin{equation}\label{Phi_nu}
\Phi_\nu=(P/Q)^\nu
\end{equation}
where $P$ and $Q$ are positive\footnote{By positive we mean that $P>0$ for any point on the $d$-torus. In contrast, if $P\geq0$ and $P=0$ for certain points, we say that $P$ is \emph{nonnegative}.} trigonometric polynomials associated with the Lagrange multipliers.
However, it is highly nontrivial to prove that the optimal dual variable lies in the interior of the feasible set. In order to overcome this difficulty, we introduced a regularization term in the dual objective function which depends solely on $P$:
\begin{equation}
\frac{\lambda}{\nu-1} \int_{\Tbb^d}\frac{1}{P^{\nu-1}}\d\m
\end{equation}
where $\lambda>0$ is a regularization parameter. Such a regularizer can of course be interpreted as a \emph{barrier function} (for $P$) since under the regularity condition 
\begin{equation}\label{old_cond_nu}
\nu\geq d/2+1, 
\end{equation}
the regularizer takes an infinite value if $P$ has a zero on the $d$-torus. Hence, the optimal $P$ is forced to be an interior point, i.e., a positive polynomial. Obviously, one can always choose $\nu$ such that the condition \eqref{old_cond_nu} is met, and indeed in this case, we showed that the optimal $Q$ is also positive so that the optimal form \eqref{Phi_nu} is true for the primal problem \eqref{primal_prob}. The price to pay for the regularization is that the $\nu$-cepstral constraints \eqref{constraint_cep} are only approximately satisfied with an error that decreases as $\lambda\to 0_+$, i.e., the rational function \eqref{Phi_nu} represents an approximate solution to \eqref{primal_prob}. It is worth noting that a process with the spectrum (\ref{Phi_nu}) can be generated by the cascade of $\nu$ identical multidimensional filters (see Fig.~\ref{fig:cascade_linear_system}) if $P/Q$ admits a spectral factor corresponding to a state space realization. Clearly, the larger $\nu$ is, the larger the complexity of the realization is. Accordingly, the key point from a practical perspective is to have the possibility to take $\nu$ as small as possible.

Having in mind this important requirement, in this letter we want to face the following problem.
\begin{problem}
 	Take the barrier function in the dual problem as 
\begin{equation}\label{regularizer_new}
\frac{\lambda}{\nu-1} \int_{\Tbb^d}\frac{1}{P^{\nu}}\d\m
\end{equation} 
and show that the corresponding approximate solution of the primal problem exists under the weaker condition $\nu\geq d /2$.
\end{problem}

\begin{example} 
	\begin{figure}[!t]
	\centering
	\includegraphics[width=0.45\textwidth]{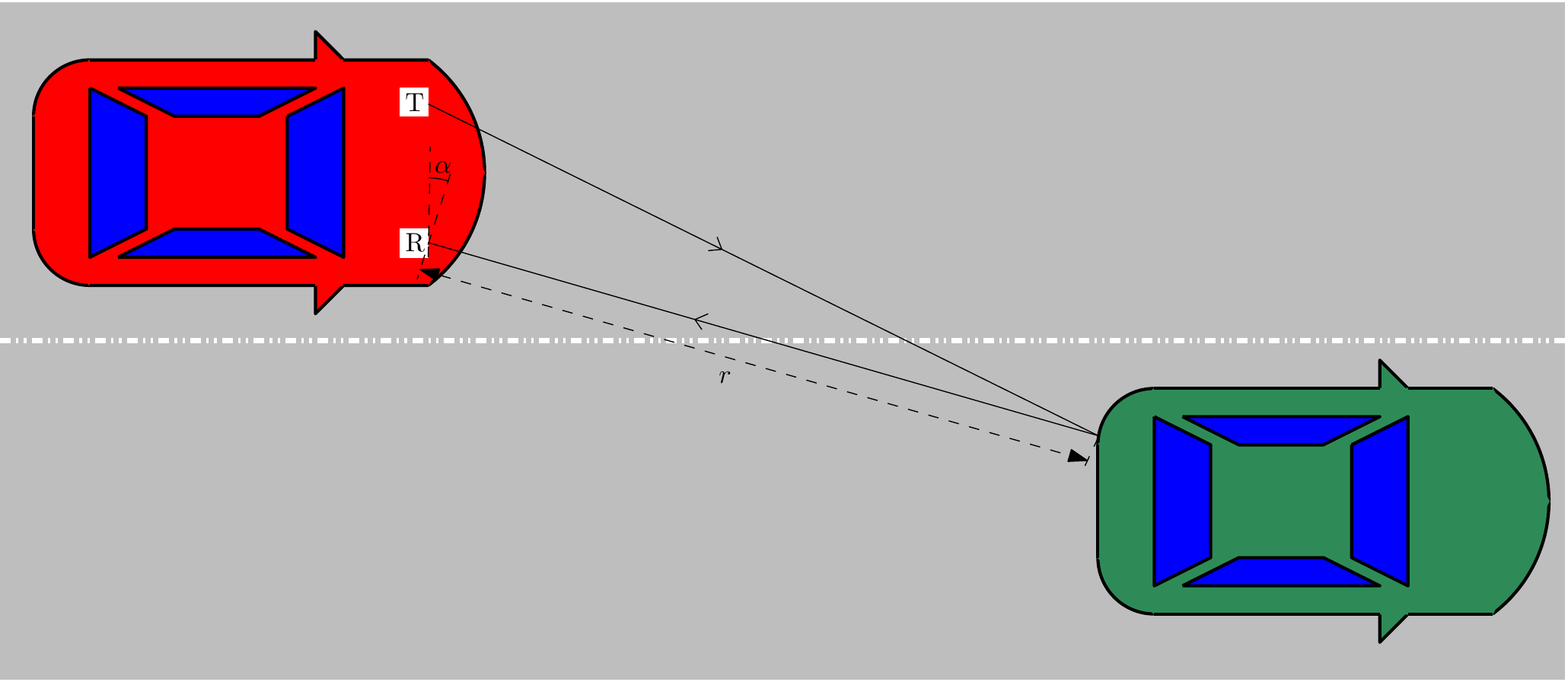}
	\caption{Target parameters estimation problem: T is source generating the pulse train signals, R is the receiver, $\alpha$ is the azimuth angle and $r$ is the range.}
	\label{fig:car}
\end{figure}
Consider an automotive radar system installed in the red car of Fig.~\ref{fig:car} that employs coherent linear frequency-modulated pulse trains signals (T) and uses a uniform linear array  of receive antennas for the measurement (R). The target (green car in Fig.~\ref{fig:car}) is identified by the range $r$, the azimuth angle $\alpha$ and the relative velocity $v$. The problem of estimating the target parameters can be formulated as a multidimensional spectral estimation problem with $d=3$ \cite{zhu2020m,engels2014target}. Under the aforementioned hypothesis on the existence of a ``realizable'' spectral factor, the estimated  model can be realized by a cascade of at least $3$ identical filters using the theory developed in \cite{Zhu-Zorzi2023cepstral}, while only $2$  using the theory that we will develop in this letter which makes more efficient the implementation of the model for simulation purposes.
\end{example}

\section{Dual problem} \label{sec:dual}

Let us recall the innocuous condition $\nu\geq 2$ which is assumed throughout this letter.
In \cite[Section~4]{Zhu-Zorzi2023cepstral}, it has been shown that the dual function of the primal problem \eqref{primal_prob} is
\begin{equation}\label{unreg_dual_func}
J_{\nu}(P,Q) := \frac{1}{\nu-1} \int_{\Tbb^d} \frac{P^\nu}{Q^{\nu-1}} \d\mu + \innerprod{\qb}{\cb} - \innerprod{\pb}{\mb},
\end{equation}
where:
\begin{itemize}
	\item $\qb=\{q_{\kb}\,:\,\kb\in\Lambda\}$ contains the real coefficients of the nonnegative trigonometric polynomial 
	\begin{equation}
	Q(e^{i\thetab}):=\sum_{\kb\in\Lambda} q_{\kb}e^{-i\innerprod{\kb}{\thetab}}
	\end{equation}
	such that $q_{-\kb}=q_{\kb}$;
	\item similarly, $\pb=\{p_{\kb}\,:\,\kb\in\Lambda_0\}$ contains the real coefficients of the nonnegative polynomial 
	\begin{equation}
	P(e^{i\thetab})=\sum_{\kb\in\Lambda} p_{\kb}e^{-i\innerprod{\kb}{\thetab}}
	\end{equation}
	with $p_{-\kb}=p_{\kb}$ and $p_{\zerob}=1$ \emph{fixed};
	\item $\innerprod{\qb}{\cb}:=\sum_{\kb\in\Lambda} q_{\kb} c_{\kb}$ inner product of two multisequences indexed in $\Lambda$, and $\innerprod{\pb}{\mb}$ is understood similarly with the index set replaced by $\Lambda_0$.
\end{itemize}
Since it is rather difficult to prove that the dual problem admits an interior-point solution, see \cite{Zhu-Zorzi2023cepstral}, we consider the regularized dual function
\begin{subequations}\label{reg_dual_func}
\begin{align}
 & J_{\nu,\lambda}(P,Q) := 
 J_{\nu}(P, Q) + \frac{\lambda}{\nu-1} \int_{\Tbb^d}\frac{1}{P^\nu}\d\mu \label{reg_dual_func_with_barrier}\\
 &\quad = \frac{1}{\nu-1} \int_{\Tbb^d} g(P(e^{i\thetab}), Q(e^{i\thetab})) \d\mu  + \innerprod{\qb}{\cb} - \innerprod{\pb}{\mb} \label{reg_dual_func_with_g} 
\end{align}
\end{subequations}
where the bivariate function is defined as
\begin{equation}\label{func_g}
g(x,y) := x^\nu/y^{\nu-1} + \lambda/x^{\nu}
\end{equation}
with $x>0,\ y>0$ and the regularization parameter $\lambda>0$ is fixed. Now let us introduce the feasible sets
\begin{equation}
		\begin{aligned}
			\Pfrak_{+} & := \left\{Q(e^{i\thetab})=\sum_{\kb\in\Lambda} q_{\kb}e^{-i\innerprod{\kb}{\thetab}} : Q>0 \ \text{on}\ \Tbb^d \right\}, \\
			\Pfrak_{+,o} & := \left\{P(e^{i\thetab})=\sum_{\kb\in\Lambda} p_{\kb}e^{-i\innerprod{\kb}{\thetab}} \in \Pfrak_{+} : p_{\zerob}=1 \right\},\nn
		\end{aligned}
	\end{equation}
so that $Q\in \Pfrak_{+}$ and $P\in \Pfrak_{+,o}$.
The domain of definition of $J_{\nu,\lambda}$ can be extended to the boundary of the feasible set $\Pfrak_{+,o}\times \Pfrak_{+}$ by excluding the zero sets of $P$ and $Q$ from the domain of integration, which does not change the values of the integrals since the zero sets have zero Lebesgue measure. Moreover, $J_{\nu,\lambda}$ may take a value of $\infty$ at some boundary points, and hence it is understood as an \emph{extended real-valued} function. 

\begin{lemma}\label{lem_strict_convex_g}
	The function $g$ in \eqref{func_g} is strictly convex in the domain $x>0,\ y>0$ $($the first quadrant$)$.
\end{lemma}
\begin{proof}
	We shall prove the claim via the derivative test. After some straightforward computations, we arrive at
	\begin{subequations}
	\begin{align}
	\nabla g(x,y) = & x^{\nu-1} y^{-\nu} \bmat \nu y \\ (1-\nu) x\emat + \bmat -\lambda\nu x^{-\nu-1} \\ 0 \emat, \label{gradient_g} \\
	\nabla^2 g(x,y) = & \nu(\nu-1) x^{\nu-2} y^{-\nu-1} \bmat y^2 & -xy \\ -xy & x^2 \emat \nn \\
	 & + \bmat \lambda\nu(\nu+1) x^{-\nu-2} & 0 \\
	 0 & 0 \emat, \label{Hessian_g}
	\end{align}
	\end{subequations}
    where $\nu(\nu-1)>0$ is a positive integer.
	It is readily observed that every diagonal element in the Hessian of $g$ is positive, and the first matrix in \eqref{Hessian_g} is positive semidefinite. After checking the determinant of the Hessian, we conclude that $\nabla^2 g(x,y)$ is positive definite in the first quadrant and the strict convexity follows.
\end{proof}

The next proposition is a direct consequence of Lemma~\ref{lem_strict_convex_g}.

\begin{proposition}\label{prop_convex}
	The regularized dual function $J_{\nu,\lambda}$ is strictly convex in the closed set $\overline{\Pfrak}_{+,o}\times\overline{\Pfrak}_+$.
\end{proposition}

\begin{proof}
	We only need to show the strict convexity of the integral term in \eqref{reg_dual_func_with_g} since the inner products are linear in $(P, Q)$. In the interior of the feasible set, namely $P\in\Pfrak_{+,o}$ and $Q\in\Pfrak_{+}$, the strict convexity of the integral term follows from that of $g$ (see Lemma~\ref{lem_strict_convex_g}), and this can be seen by a pointwise argument on the integrand (see e.g., Proposition~5.3 in \cite{zhu2020m}). The same reasoning works if $P$ or $Q$ is on the boundary of the respective feasible set because, once the zero sets of $P$ and $Q$ are excluded from the domain of integration, the function $g(P(e^{i\thetab}), Q(e^{i\thetab}))$ is well defined and the proof holds verbatim.	
\end{proof}

We conclude that the regularized dual optimization can be formulated as
\begin{equation}\label{reg_dual_prob}
\min\ J_{\nu, \lambda}(P, Q)\quad \text{s.t.}\ P\in \overline\Pfrak_{+,o},\ Q\in \overline\Pfrak_{+}.
\end{equation}

\section{A unique interior-point solution under the condition $\nu\geq {d}/{2}$}\label{sec:weaker_cond}

In this section we shall see how the condition $\nu\geq d/2$ guarantees that the regularized dual problem (\ref{reg_dual_prob}) has an interior-point solution. By Proposition~\ref{prop_convex}, we know that a minimizer of $J_{\nu, \lambda}$ in $\overline{\Pfrak}_{+,o}\times\overline{\Pfrak}_+$ is unique \emph{provided that it exists}. We shall first establish such existence. Then, under the regularity condition $\nu\geq d/2$ which is weaker than \eqref{old_cond_nu}, we aim to exclude the possibility that a minimum of $J_{\nu,\lambda}$ may fall on the boundary of the feasible set, which is a weaker version of Lemma~5.8 in \cite{Zhu-Zorzi2023cepstral}, using a Byrnes--Gusev--Lindquist-type argument that first appeared in \cite{BGL-98} and subsequently in e.g., \cite{RKL-16multidimensional}. 

\subsection{Existence of a minimizer in $\overline{\Pfrak}_{+,o}\times\overline{\Pfrak}_+$}\label{subsec:exist}

The existence of a solution to \eqref{reg_dual_prob} can be shown via reasonings similar to the ones in \cite{Zhu-Zorzi2023cepstral}: such existence depends on the following feasibility assumption.

\begin{assumption}[Feasibility]\label{assump_feasible}
	The given covariances $\{c_\kb\}_{\kb\in\Lambda}$ admit an integral representation
	\begin{equation}
	c_\kb = \int_{\Tbb^d} e^{i\innerprod{\kb}{\thetab}} \Phi_0 \d\m\quad \forall \ \kb\in\Lambda, 
	\end{equation}
	where $\Phi_0$ is a nonnegative function on $\Tbb^d$ and is positive on some open ball $B_1\subset\Tbb^d$.
\end{assumption}

In order to prove our existence result, we need the following lemmas.
\begin{lemma}\label{lem_lower_semicont}
	The unregularized dual function $J_\nu$ in \eqref{unreg_dual_func} and the regularized version $J_{\nu,\lambda}$ in \eqref{reg_dual_func} are lower-semicontinuous on $\overline{\Pfrak}_{+,o}\times\overline{\Pfrak}_+$. In particular, they are both continuous on $\Pfrak_{+,o}\times\Pfrak_+$.
\end{lemma}
\begin{proof}
	The proof is similar to that of Lemma~5.5 in \cite{Zhu-Zorzi2023cepstral}. The only difference here is that the power of $P$ in the regularizer in \eqref{reg_dual_func_with_barrier} is $\nu$ instead of $\nu-1$ in \cite{Zhu-Zorzi2023cepstral}.
\end{proof}

\begin{lemma}\label{lem_unbounded}
	Suppose that Assumption~\ref{assump_feasible} holds. If a sequence $\{(P_j,Q_j)\}_{j\geq1}\subset\overline{\Pfrak}_{+,o}\times \overline{\Pfrak}_+$ is such that $\|(P_j,Q_j)\|\to\infty$ as $j\to\infty$, then $J_{\nu,\lambda}(P_j,Q_j)\to\infty$.		
\end{lemma}
\begin{proof}
	See the proof of Lemma~5.6 in \cite{Zhu-Zorzi2023cepstral} which uses Assumption~\ref{assump_feasible}.
\end{proof}

\begin{proposition}\label{prop_exist}
	Under Assumption~\ref{assump_feasible}, the regularized dual optimization problem \eqref{reg_dual_prob} admits a solution.
\end{proposition}
\begin{proof}
	Take a sufficiently large real number $\beta$, so that the  sublevel set of the regularized dual function
	\begin{equation}\label{sublevel_set_reg}
		J_{\nu,\lambda}^{-1}(-\infty,\beta] := \{ (P,Q)\in \overline{\Pfrak}_{+,o}\times \overline{\Pfrak}_+ \,:\, J_{\nu,\lambda}(P,Q)\leq\beta \}.
	\end{equation}
    is not empty. Then Lemma~\ref{lem_lower_semicont} implies that the sublevel set is closed; Lemma~\ref{lem_unbounded} implies that the sublevel set is bounded.
	Obviously the polynomial pair $(P,Q)$, parametrized by their coefficients, belongs to a finite-dimensional vector space. It follows that the sublevel set is compact.
	Given the lower-semicontinuity of the objective function $J_{\nu,\lambda}$ (see Lemma~\ref{lem_lower_semicont}), a minimizer exists in $J_{\nu,\lambda}^{-1}(-\infty,\beta]$ by the extreme value theorem of Weierstrass.
\end{proof}

\subsection{Non-optimality of boundary points given $\nu\geq d/2$}\label{subsec:boundary}

In this subsection we prove that the optimal solution of (\ref{reg_dual_prob}) cannot belong to the boundary of $\overline{\Pfrak}_{+,o}\times \overline{\Pfrak}_{+}$ using arguments which conceptually differ from the ones used for the case $\nu\geq d/2+1$ in \cite{Zhu-Zorzi2023cepstral}. Notice first that if $P\in\partial\Pfrak_{+,o}$ where $\partial$ denotes the boundary of a set, then the regularization term \eqref{regularizer_new}
employed in \eqref{reg_dual_func_with_barrier} takes a value of $\infty$ under the condition $\nu\geq d/2$, see Proposition~A.4 in \cite{zhu2020m}. The other term, namely the unregularized function $J_\nu(P,Q)$ whose expression is given in \eqref{unreg_dual_func}, is bounded from below under Assumption~\ref{assump_feasible}, see the proof of \cite[Lemma~5.6]{Zhu-Zorzi2023cepstral}. Therefore, in this case we have $J_{\nu,\lambda}(P,Q)=\infty$ which is certainly not a minimum. Consequently, an optimal $(P, Q)$ must have $P\in\Pfrak_{+,o}$.

Next, we work on the case of $(P,Q)\in\Pfrak_{+,o}\times\partial\Pfrak_{+}$. It is still possible that the integral term in \eqref{unreg_dual_func} diverges so that $J_{\nu,\lambda}(P, Q)=\infty$, and such a point is obviously not a minimizer. Therefore, we only need to consider points $(P, Q)\in\Pfrak_{+,o}\times\partial\Pfrak_{+}$ such that the function value $J_{\nu,\lambda}(P,Q)$ is \emph{finite}. Construct the real-valued function 
\begin{equation}\label{func_f}
	f(t) := J_{\nu,\lambda}(P, Q+t\oneb)
\end{equation}
defined for $t\geq 0$ where $\oneb$ denotes the constant polynomial taking value one.
One can show without difficulty that $f$ ``inherits'' the strict convexity from $J_{\nu,\lambda}$ (see Proposition~\ref{prop_convex}). It is worth noting that the continuity of $f(t)$ can be extended to $t=0_+$ which corresponds to the boundary point $(P, Q)$, in contrast with the general result that $J_{\nu,\lambda}$ is only lower-semicontinuous on the boundary of the feasible set (see Lemma~\ref{lem_lower_semicont}). These properties are established in the next statement.

\begin{lemma}\label{lem_bdry_cont}
	The function $f(t)$ in \eqref{func_f}, where $(P,Q)$ is an arbitrarily fixed point in $\Pfrak_{+,o}\times\partial\Pfrak_{+}$ such that $J_{\nu,\lambda}(P,Q)<\infty$, is strictly convex and continuous in $[0, \infty)$.
\end{lemma}
\begin{proof}
	The strict convexity of $f$ follows directly from Proposition~\ref{prop_convex}, and it remains to show the continuity.	
	For $t>0$, the argument $(P, Q+t\oneb)$ belongs to the interior $\Pfrak_{+,o}\times \Pfrak_{+}$, and the continuity of $f$ follows from that of $J_{\nu,\lambda}$ (see Lemma~\ref{lem_lower_semicont}).
	Hence, we are only concerned with the right continuity of $f$ at $t=0$. 
	By the lower-semicontinuity of $J_{\nu,\lambda}$ at $(P,Q)$
	(Lemma~\ref{lem_lower_semicont} again), for any $\varepsilon>0$, there exists $\delta>0$ such that whenever $(P_1,Q_1)\in \overline{\Pfrak}_{+,o}\times\overline{\Pfrak}_+$ satisfies $\|(P_1, Q_1) - (P, Q)\|<\delta$, we have $J_{\nu,\lambda}(P_1,Q_1)>J_{\nu,\lambda}(P,Q)-\varepsilon$. We can always choose $t$ sufficiently small to make the argument on the right-hand side of \eqref{func_f} sufficiently close to $(P,Q)$, which leads to the inequality $f(t)> f(0)-\varepsilon$. At the same time, by strict convexity we have $f(t)<tf(1)+(1-t)f(0) = f(0)+t[f(1)-f(0)]$.
	The last term, namely $t[f(1)-f(0)]$ can be made smaller than $\varepsilon$ (in absolute value) for $t$ sufficiently small. Therefore, we reach the inequality $|f(t)-f(0)|<\varepsilon$ which proves the continuity at $t=0_+$ since $\varepsilon>0$ is arbitrarily chosen.
\end{proof}


\begin{proposition}\label{prop_boundary}
	Fix any point $(P,Q)\in\Pfrak_{+,o}\times\partial\Pfrak_{+}$ such that $J_{\nu,\lambda}(P,Q)<\infty$. Then the function $f(t)$ in \eqref{func_f} has a derivative	$f'(t)\to -\infty$ as $t\to 0_+$. Therefore, $(P,Q)$ cannot be a minimizer of $J_{\nu,\lambda}$.
\end{proposition}
\begin{proof}
	After some computations, we have for $t>0$ that
	\begin{equation}\label{func_f_deriva}
		f'(t) = c_{\zerob} - \int_{\Tbb^d} \left[ \frac{P}{Q+t\oneb}\right]^\nu \d\mu.
	\end{equation}
    Since $P\in\Pfrak_{+,o}$, we know $P_{\min}:=\min_{\thetab\in\Tbb^d} P(e^{i\thetab})>0$.
	Then the following relation
	\begin{equation}
		\int_{\Tbb^d} \left[ \frac{P}{Q+t\oneb}\right]^\nu \d\mu \geq P_{\min}^\nu \int_{\Tbb^d}  \frac{1}{(Q+t\oneb)^\nu} \d\mu \to\infty
	\end{equation}
    holds as $t\to 0_+$ by Lebesgue's monotone convergence theorem \cite[p.~21]{rudin1987real} and Proposition~A.4 in \cite{zhu2020m}.
    This shows that $f'(t)$ in \eqref{func_f_deriva} tends to $-\infty$ as $t$ goes to zero from the right. Taking Lemma~\ref{lem_bdry_cont} into account, we conclude that $0$ is not a local minimizer of $f(t)$. Consequently, any boundary point $(P,Q)\in\Pfrak_{+,o}\times\partial\Pfrak_{+}$ is not a minimizer of $J_{\nu,\lambda}$, because taking an arbitrarily small step along the direction $(\zerob,\oneb)$, which points towards the interior $\Pfrak_{+,o}\times\Pfrak_{+}$, will result in a decrease of the objective function value.
\end{proof}

We summarize what we have got so far in the next theorem.
\begin{theorem}\label{thm_main}
	Under Assumption~\ref{assump_feasible} and the condition $\nu>d/2$, the optimization problem \eqref{reg_dual_prob} admits a unique interior-point solution $(\hat{P}, \hat{Q})\in\Pfrak_{+,o}\times \Pfrak_{+}$ such that
	\begin{subequations}\label{stationarity_cond}
		\begin{align}
			c_\kb & = \int_{\Tbb^d} e^{i\innerprod{\kb}{\thetab}} ({\hat P}/{\hat Q})^{\nu}\d\m \; \;  \forall\kb\in\Lambda, \label{part_Q} \\
			m_\kb & = \int_{\Tbb^d} e^{i\innerprod{\kb}{\thetab}} \frac{\nu}{\nu-1}\left[({\hat P}/{\hat Q})^{\nu-1}-{\lambda}/{\hat P^{\nu+1}} \right] \d\m  \; \;  \forall\kb\in\Lambda_0. \label{part_P}
		\end{align}
	\end{subequations}
    In plain words, the spectral density $\hat{\Phi}_\nu=(\hat{P}/\hat{Q})^\nu$ achieves covariance matching and approximate $\nu$-cepstral matching with errors
    \begin{equation*}
    	\varepsilon_{\kb} = \frac{\lambda\nu}{\nu-1}\int_{\Tbb^d} e^{i\innerprod{\kb}{\thetab}} \frac{1}{\hat P^{\nu+1}} \d\m,\ \kb\in\Lambda_0.
    \end{equation*}
\end{theorem}

\begin{proof}
	The existence of a solution is guaranteed by Proposition~\ref{prop_exist} and the uniqueness by Proposition~\ref{prop_convex}. Moreover, given the reasoning at the beginning of this subsection and Proposition~\ref{prop_boundary}, the optimal $(\hat{P}, \hat{Q})$ must be an interior point, i.e., both polynomials are positive. As a consequence, it must satisfy the stationary-point equation $\nabla J_{\nu,\lambda}(P, Q) = \zerob$, which is equivalent to the conditions in \eqref{stationarity_cond}. 
	Indeed, this point can be seen by setting the first differential of the regularized dual function
	\begin{equation}\label{first_diff_reg_dual_func}
	\begin{aligned}
	 & \delta J_{\nu,\lambda}(P, Q; \delta P, \delta Q) = \innerprod{\cb}{\delta\qb} - \int_{\Tbb^d} \delta Q (P/Q)^\nu \d\mu \\
	 & \quad - \innerprod{\mb}{\delta\pb} + \frac{\nu}{\nu-1} \int_{\Tbb^d} \delta P\left[ (P/Q)^{\nu-1} - \lambda/P^{\nu+1}\right] \d\m
	\end{aligned}
	\end{equation}
    equal to zero for any direction $(\delta P, \delta Q)$.
\end{proof}

\begin{remark}\label{remark1}
It is important to observe that the continuous dependence of the solution $(\hat P, \hat Q)$ to \eqref{reg_dual_prob} on the covariance and $\nu$-cepstral data $(\cb, \mb)$ can be established similarly to \cite[Sec.~6]{Zhu-Zorzi2023cepstral}, so that the optimization problem \eqref{reg_dual_prob} is in fact well-posed in the sense of Hadamard.
\end{remark}

\section{Numerical simulations}\label{sec:nume}


In this section we present some numerical experiments in which the problem \eqref{primal_prob} is used to reconstruct the spectrum of a $3$-d stationary random field $y(t_1, t_2, t_3)$ starting from a finite set of its covariance lags and $\nu$-cepstral coefficients. We assume that the underlying process $y$ is described as the output of a cascade linear shaping filter $\tilde{W}(z_1, z_2, z_3)=W^\nu(z_1,z_2,z_3)$ with a white noise input $e(t_1, t_2, t_3)$, see Fig.~\ref{fig:cascade_linear_system} with $\nu=2$ and $d=3$. We shall also assume that the transfer function $W$ has the structure consistent with our optimal form \eqref{Phi_nu} for the spectrum. More precisely, we consider the class of rational models such that the numerator and denominator polynomials both have degree one:
\begin{equation}\label{transfer_func_W}
	W(\zb) = \left[ \frac{b(\zb)}{a(\zb)} \right]^\nu = \left[ \frac{b_0 - b_1z_1^{-1} - b_2z_2^{-1} - b_3z_3^{-1}}{a_0 - a_1z_1^{-1} - a_2z_2^{-1} - a_3z_3^{-1}} \right]^\nu
\end{equation}
where $\zb$ represents $(z_1,z_2,z_3)$ for short. Obviously, the polynomials $a(\zb)$ and $b(\zb)$ are described by the respective vectors 
$\ab=[a_0,\dots,a_3]$ and $\bb=[b_0,\dots,b_3]$ of coefficients.
If the white noise input $e$ has unit variance, then the spectral density of the output process $y$ is $(P/Q)^\nu$ where
$P(e^{i\thetab}) = |b(e^{i\thetab})|^2$ and $Q(e^{i\thetab}) = |a(e^{i\thetab})|^2$ 
are two nonnegative polynomials in $\overline{\Pfrak}_+$. Throughout this section, we set the parameter $\nu=2$ to meet the condition $\nu\geq {d}/{2}=3/2$.
In such scenario it is worth stressing that the theory developed in  \cite{Zhu-Zorzi2023cepstral} does not work. Indeed that theory requires $\nu\geq 3$ to guarantee the existence of an approximate solution to the primal problem \eqref{primal_prob}.

In what follows we consider two models of the form \eqref{transfer_func_W}. We take two sets of \emph{real} parameters $(\ab_j, \bb_j)$ with $j=1, 2$ such that
\begin{equation}
	\begin{aligned}
	\ab_1 = \ab_2 =\ab & = [1, 0.3, 0.3, 0.3], \\
	\tilde{\bb}_1 = [1, -0.2, -0.3, -0.4], \quad & \tilde{\bb}_2 = [1, -0.2, -0.3, -0.5],
	\end{aligned}
\end{equation}
and $\bb_j=\tilde{\bb}_j/\|\tilde{\bb}_j\|$. The last operation of normalization gives $\|\bb\|^2=1$ which is equivalent to $p_\zerob=1$ for the numerator polynomial $P$. The first model, hereafter called ``zeroless model'', corresponds to the process with spectrum $\Phi_1=(P_1/Q)^2$ where the polynomials $Q=|a|^2$ and $P_1=|b_1|^2$ are positive on $\Tbb^3$; the second model, hereafter called ``model with a spectral zero'', corresponds to the process with spectrum $\Phi_2=(P_2/Q)^2$ where  $P_2=|b_2|^2$ (and also the spectrum) has a zero at the frequency vector $(\pi,\pi,\pi)$. The index set $\Lambda$ is identified as $\Lambda = \Lambda_+ - \Lambda_+$ with
\begin{equation}
\Lambda_+ := \{\, (0,0,0),\, (1,0,0),\, (0,1,0),\, (0,0,1) \,\}.
\end{equation}
Here the set difference is understood as $A-B:=\{x-y : x\in A,\ y\in B\}$. 
Since $q_{-\kb}=q_{\kb}$ and $p_{-\kb}=p_{\kb}$, the total number of variables is $13$.

By means of the discrete formulation described in Section~7 of \cite{Zhu-Zorzi2023cepstral}, we consider the following procedure to test the ability to reconstruct the spectra through the solution of \eqref{primal_prob} for the previous two models:
\begin{enumerate}
	\item Fix $N=20$ and discretize $\Tbb^3$ into $N^3$ regular grid points by gridding the interval $[0, 2\pi]$ into $N$ equidistant points in each dimension. 
	\item Compute the ``true'' covariances $\{c_\kb : \kb\in\Lambda\}$ and the $\nu$-cepstral coefficients $\{m_\kb : \kb\in\Lambda_0\}$ of $\Phi=|W|^2$ via \eqref{constraint_cov} and \eqref{constraint_cep} where $\d\m$ is replaced by a discrete measure with equal mass $1/N^3$ on the grid points. 
	\item Solve the discrete version of the regularized dual problem \eqref{reg_dual_prob} using $(\cb, \mb)$ computed above and  $\lambda>0$ chosen sufficiently small. 
	\item Let $(\hat{\pb},\hat{\qb})$ be the optimal solution to \eqref{reg_dual_prob}
	and $(\pb,\qb)$ be the polynomial coefficients  corresponding to the true spectrum $\Phi=|W|^2$. Finally, evaluate the reconstruction error $\|(\hat{\pb},\hat{\qb}) - (\pb,\qb)\|$.
\end{enumerate}

In Step 3, the optimization problem is solved using Newton's method. Some computational details can be found in \cite{Zhu-Zorzi2023cepstral}, and suitable modifications on the gradient and Hessian of the objective function are needed for the current problem \eqref{reg_dual_prob}. 

\begin{figure}[!t]
	\centering
	\includegraphics[width=4.1cm]{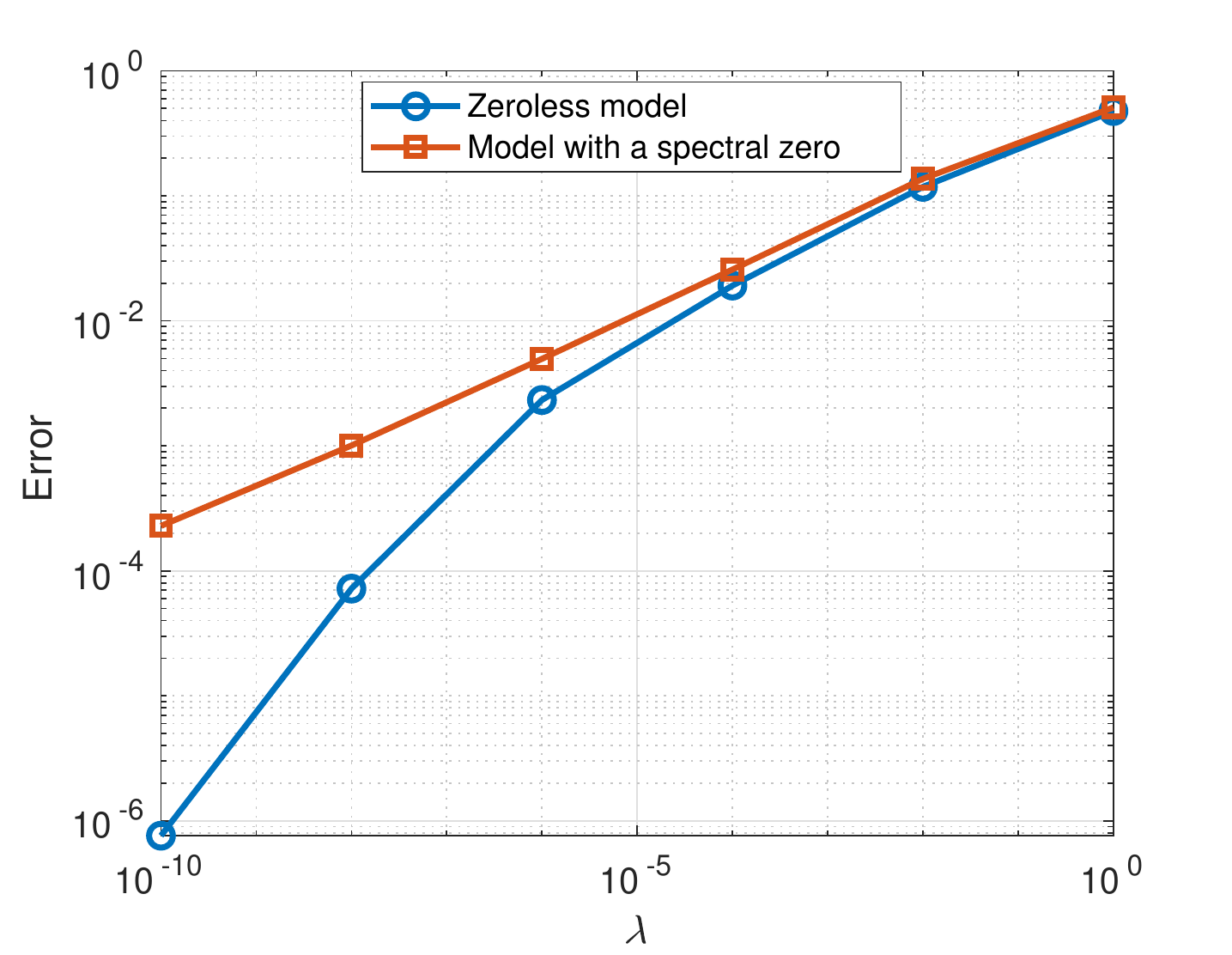}
	\includegraphics[width=4.4cm]{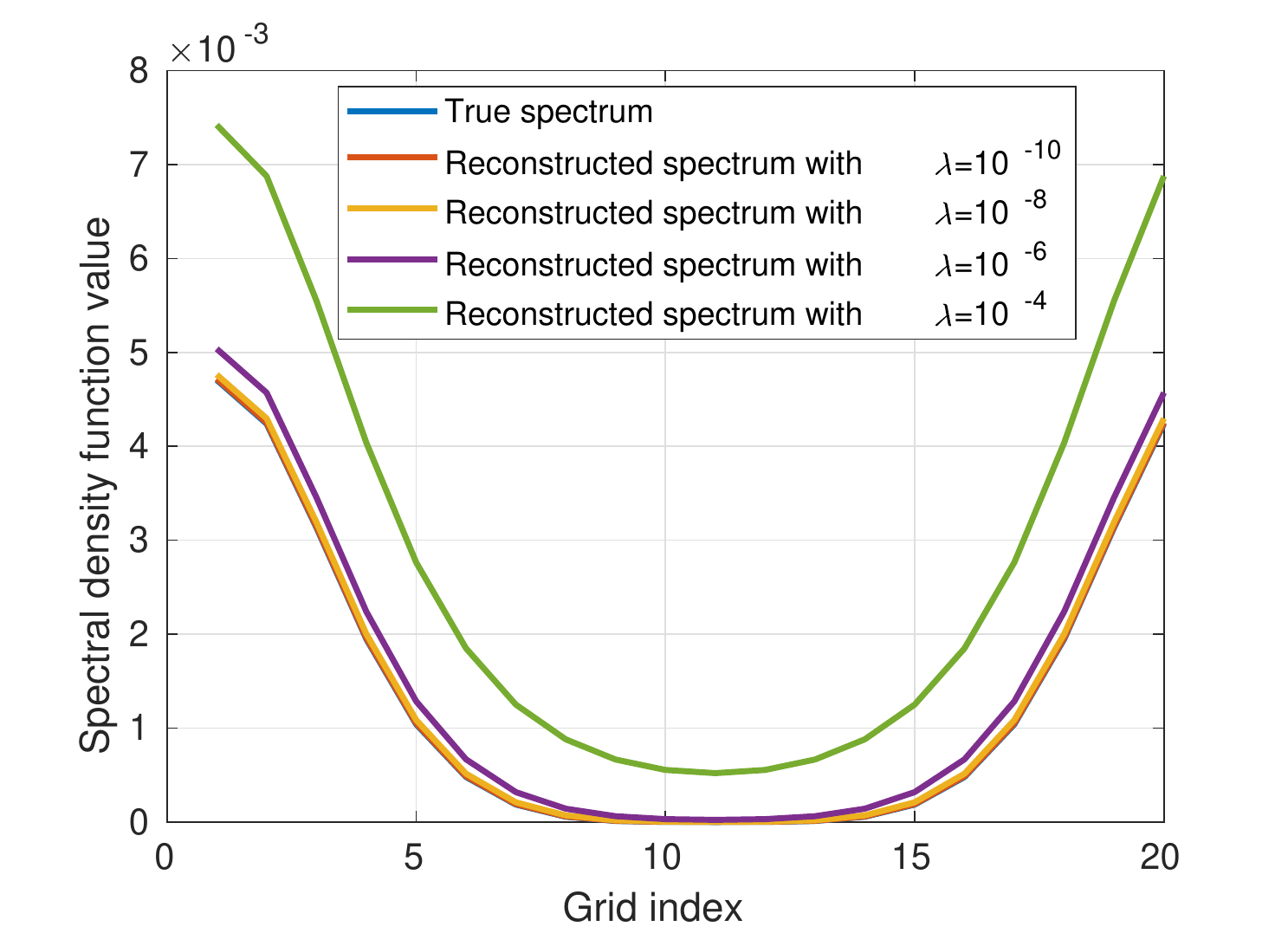}
	\caption{\emph{Left:} Error of spectrum reconstruction versus the regularization parameter $\lambda = 10^{-n}$, $n=0, 2, 4, 6, 8, 10$, where both axes are in the logarithmic scale. \emph{Right:} The true spectrum with a zero and the reconstructed spectra with choices of $\lambda=10^{-10}, 10^{-8}, 10^{-6}, 10^{-4}$ at the cross section $[\,\cdot \;11\;11\,]$, i.e., $\hat{\Phi}_{\nu,\lambda}(e^{i\thetab})$ with $\thetab=2\pi\times[(k-1)/20, 11/20, 11/20]$ for grid indices $k=1,\ldots, 20$. Note that the reconstructed spectra with $\lambda=10^{-10}$ and $10^{-8}$ (red and yellow lines) and the true spectrum (blue line) almost overlap each other.}
	\label{fig:secx_vs_lambda}
\end{figure}

The left panel of Fig.~\ref{fig:secx_vs_lambda} shows the reconstruction errors for the two models with different values of the regularization parameter $\lambda= 10^{-n}$, $n=0, 2, 4, 6, 8, 10$. It is readily observed in both cases that the errors decrease monotonically as $\lambda\to 0$.  
In view of Remark~\ref{remark1}, therefore, if we have a dataset generated by a $d$-dimensional process and consider the approximate solution to the problem (\ref{primal_prob}) with $\nu=\lceil {d/2}\rceil$ where $\cb$ and $\mb$ are replaced by their sample estimators (computed from the data), then the resulting spectral estimator is characterized by a small estimation error provided that $\lambda$ is chosen sufficiently small and the underlying process has a spectrum of the form (\ref{Phi_nu}) with $\nu=\lceil {d/2}\rceil$ which agrees with our parameter specification.

In the right panel of Fig.~\ref{fig:secx_vs_lambda}, we compare the spectral density of the model with a spectral zero and the reconstructed spectra for values of $\lambda=10^{-10}, 10^{-8}, 10^{-6}, 10^{-4}$ (corresponding to the orange line in the left panel) along a cross section $[\,\cdot \;11\;11\,]$ of the regular grid for $\Tbb^3$, where the true spectral zero is located at the index\footnote{In Matlab, the array indices start with $1$.} $[11, 11, 11]$.
Notice that the other two cross sections $[\,11\;\cdot \;11\,]$ and $[\,11\;11\;\cdot \,]$ of the spectral densities are not shown because they are visually similar to this figure. We conclude that as $\lambda$ goes to zero, the true nonnegative spectrum is well approximated by a positive spectral density with smaller and smaller errors, which is consistent with the left panel.

\section{Conclusions}\label{sec:concl}
We have considered a multidimensional $\nu$-moment problem which searches a spectral density maximizing the $\nu$-entropy and matching a finite set of covariance lags and $\nu$-cepstral coefficients. A previous work showed that it is possible to guarantee the existence of a rational approximate solution only in the case $\nu\geq  d/2+1$ where $d$ is the dimension of the random field described by the multidimensional spectrum. Motivated by the fact that $\nu$ should be chosen as small as possible, we have proposed a different regularization technique which ensures the existence of a rational approximate solution to the $\nu$-moment problem under the weaker condition $\nu\geq  d/2$.

%


%
%
%
%

\bibliographystyle{IEEEtran}
\bibliography{references}   

\end{document}